\documentclass[11pt, english]{smfart}
\usepackage{amsfonts, amssymb, amsmath, latexsym}
\usepackage{mathptmx}
\usepackage[english]{babel}
\usepackage{mysectionsen}
\usepackage[latin1]{inputenc}
\usepackage[T1]{fontenc}
\usepackage{textcomp}

\usepackage{graphicx}
\usepackage[dvips]{color}
\input xy
\xyoption{all}

\newtheorem{Example}[Thm]{Example}

\newtheorem{Remark}[Thm]{Remark}
\newtheorem{Remarks}[Thm]{Remarks}

\newcounter{thm}

\newtheorem{main_theorem}[thm]{Theorem}

\begin{document}

\thispagestyle{empty}

\begin{center}

{\Large \bfseries{\textsc{An intrinsic characterization of Bruhat-Tits buildings inside analytic groups}}}

\end{center}

\vspace{0.3cm}

\begin{center}
\textsc{Bertrand R\'emy, Amaury Thuillier and Annette Werner}
\end{center}

\vspace{3cm}
\hrule
\vspace{0,5cm}

{\small
\noindent
{\bf Abstract:}~Given a semisimple group over a complete non-Archimedean field, it is well known that techniques from non-Archimedean analytic geometry provide an embedding of the corresponding Bruhat-Tits builidng into the analytic space associated to the group; by composing the embedding with maps to suitable analytic proper spaces, this eventually leads to various compactifications of the building. 
In the present paper, we give an intrinsic characterization of this embedding. 

\vspace{0,1cm}


\vspace{0,5cm}
\hrule
\vspace{0,5cm}

\noindent
{\bf AMS classification (2000):}
20E42,
51E24,
14L15,
14G22.
}

\vspace{0,5cm}
\hrule

\vspace{3cm}

{\it Dedicated to Gopal Prasad with admiration}

\vspace{3cm}

\tableofcontents

\newpage

\section*{Introduction}
\label{s - intro}

Bruhat-Tits theory (\cite{BT1a}, \cite{BT1b}) deals with reductive groups over valued fields. 
Given such a group ${\rm G}$ over a complete non-Archimedean field $k$, the theory provides a certain cell complex $\mathcal{B}({\rm G},k)$, called a \emph{Euclidean} or \emph{affine building}, on which the group of rational points ${\rm G}(k)$ acts in a very balanced way (see \cite[Introduction]{BT1b} for more precise conditions on the ground field). 
The space $\mathcal{B}({\rm G},k)$ is large enough with respect to the group ${\rm G}(k)$ in the sense that the group action is proper ({\it i.e.} cell stabilizers are compact subgroups), and the group ${\rm G}(k)$ is large enough with respect to the space $\mathcal{B}({\rm G},k)$  in the sense that the group action is \emph{strongly transitive}~({\it i.e.} the action is transitive on the inclusions of a maximal cell -- called an \emph{alcove} -- in a maximal flat subspace -- called an \emph{apartment}). 

\smallskip 

The latter property implicitly uses the fact that Euclidean buildings are complete metric spaces with non-positive curvature properties, which makes them relevant to geometric group theory \cite{BriHae}: this is the geometric side of Bruhat-Tits theory. 
The second aspect of the theory is of arithmetic nature. 
Namely, in order to perform the construction of the buildings, F.~Bruhat and J.~Tits define and study integral models of the group ${\rm G}$: {\it a posteriori}, each cell in the building -- called a \emph{facet} -- gives rise to such a model whose group of integral points is, possibly up to finite index, the stabilizer of the facet. 
These integral models turned out to be very useful in various mathematical fields, such as representation theory (\cite{MoyPrasad1}, \cite{MoyPrasad2}), so that the viewpoint of integral models can be chosen as main viewpoint for Bruhat-Tits theory, as shown by the upcoming book \cite{KaPra}. 

\smallskip 

Bruhat-Tits buildings are singular spaces (from the viewpoint of differential geometry for instance) since they are products of simplicial complexes. 
Therefore it was a big surprise when V.~Berkovich noticed that they could be related to non-Archimedean analytic geometry from his perspective \cite[Chapter 5]{Ber1}. 
We developed this approach and this led us to the construction of families of compactifications by analogy with those associated to Riemannian symmetric spaces (\cite{RTW1}, \cite{RTW2}, \cite{RTW4}, \cite{Chanfi}). 
For a semisimple group ${\rm G}$ over a complete non-Archimedean local field $k$, these compactifications are defined to be the closures of embedding maps from $\mathcal{B}({\rm G},k)$ to analytifications of proper ${\rm G}$-schemes (flag varieties ${\rm G}/{\rm P}$ or the wonderful compactification $\overline {\rm G}$ associated to ${\rm G}$). 
More precisely, the required maps are obtained by constructing first a canonical ${\rm G}(k)$-equivariant embedding $\vartheta : \mathcal{B}({\rm G},k) \rightarrow {\rm G}^{\rm an}$, where ${\rm G}(k)$ acts on the analytification ${\rm G}^{\rm an}$ by conjugation, and then by composing it with the analytication of the maps ${\rm G} \to {\rm G}/{\rm P}$ (where ${\rm P}$ is a suitable parabolic subgroup) or of the de Concini-Procesi map ${\rm G} \hookrightarrow \overline {\rm G}$. 

\smallskip 

In this paper, we go back to the canonical ${\rm G}(k)$-equivariant embedding $\vartheta : \mathcal{B}({\rm G},k) \rightarrow {\rm G}^{\rm an}$ itself and address the following basic problem about the relationship between Bruhat-Tits and Berkovich theories: {\it give an intrinsic description of the building $\mathcal{B}({\rm G},k)$ as a subspace of the analytic group ${\rm G}^{\rm an}$}. 
Our answer is the content of Theorem \ref{Thm-image}; it is also rephrased in Theorem \ref{characterization2}. 
It is intrinsic in the sense that the information we use from a point $z \in {\rm G}^{\rm an}$ to decide whether it belongs to the image of $\vartheta$ or not, is its \emph{holomorphic envelope} in ${\rm G}^{\rm an}$ defined as ${\rm G}(z) = \{y \in {\rm G}^{\rm an} \ ; \ \forall f \in \mathcal{O}({\rm G}), \ |f(y)| \leqslant |f(z)|\}$. 

\medskip 

\begin{main_theorem}
\label{introduction-theorem} 
The image of the canonical embedding $\vartheta : \mathcal{B}({\rm G},k) \hookrightarrow {\rm G}^{\rm an}$ is the subset of points $z \in {\rm G}^{\rm an}$ whose holomorphic envelope ${\rm G}(z)$ is a $k$-affinoid subgroup potentially of Chevalley type for which there exists a maximal torus ${\rm T}$ of ${\rm G}$ containing a maximal split torus and such that ${\rm G}(z) \cap {\rm T}^{\rm an}$ is the maximal affinoid subgroup ${\rm T}^1$ of ${\rm T}^{\rm an}$.
\end{main_theorem}

Being potentially of Chevalley type is our first (natural) necessary condition; it is explained in Section \ref{ss-Pot_Chevalley}. 
In Section \ref{ss-Galois}, we show why this condition cannot be sufficient: this is related to the fact, well known when performing descent of the ground field in Bruhat-Tits theory, that Galois fixed points sets are in general bigger than the expected rational buildings. 
The second condition deals with tori, or equivalently with apartments in the building: the existence of a torus intersecting nicely with ${\rm G}(z)$ corresponds to the fact that the point $z$ does not belong to the undesirable part of Galois fixed point sets. 
In Section \ref{ss-Recovering apartments} we explain that the combination of these two conditions is a sufficient one in order to belong to the image of $\vartheta$. 
In Section \ref{ss-Reformulation norms}, we provide a reformulation of the above theorem according to which the building of ${\rm G}(k)$ has a natural realization as a space of multiplicative $k$-norms on the coordinate ring $\mathcal{O}({\rm G})$ of ${\rm G}$. 

\medskip 

In what follows, we denote by $k$ a field which is complete with respect to a given non-Archimedean absolute value and we denote by ${\rm G}$ a semisimple group defined over $k$. 
We use the functoriality of the building $\mathcal{B}({\rm G},k)$ with respect to non-Archimedean extensions of $k$; it is satisfied at least when ${\rm G}$ is split or when $k$ is discretely valued with perfect residue field (see \cite[1.3.3]{RTW1}). 
Each of the latter two conditions also allow us to use the following fact: for any maximal split torus ${\rm S}$ in ${\rm G}$ there are a maximal torus ${\rm T}$ containing ${\rm S}$ and a finite Galois extension $k'/k$ splitting ${\rm T}$ such that, for the map given by the previous functoriality, the apartment of ${\rm T} \otimes_k k'$ contains the apartment of ${\rm S}$ (this is clear when ${\rm G}$ is split and it follows from \cite[Corollaire 5.1.12]{BT1b} in the second case). 
At last, we denote by ${\rm G}^{\rm an}$ the analytic space attached to ${\rm G}$: it is a $k$-analytic space in the sense of Berkovich whose underlying topological space consists of the multiplicative seminorms on the coordinate ring $\mathcal{O}({\rm G})$ inducing the initial absolute value on $k$. 

\medskip

{\bf Acknowledgement: }The third author was partially supported by the Deutsche Forschungsgemeinschaft (DFG, German Research Foundation) TRR 326 \textit{Geometry and Arithmetic of Uniformized Structures}, project number 444845124.

\section{Affinoid subgroups potentially of Chevalley type}
\label{ss-Pot_Chevalley}

As already mentioned in the introduction, recall \cite[Theorem 2.1]{RTW1} that to any point $x$ of
$\mathcal{B}({\rm G},k)$ can be attached a $k$-affinoid subgroup ${\rm G}_x$ of ${\rm
  G}^{\rm  an}$ satisfying the following condition: for any non-Archimedean
extension $K/k$, the subgroup ${\rm G}_x(K)$ of ${\rm
  G}(K)$ is the stabilizer of $x$ seen
in the building $\mathcal{B}({\rm G},K)$. By definition, the point
$\vartheta(x)$ is the unique element of the Shilov boundary of ${\rm
  G}_x$, {\it i.e.} the only point of ${\rm G}_x$ such
that $|f(y)| \leqslant |f(\vartheta(x))|$ for any $y \in {\rm G}_x$
and any $f \in \mathcal{O}({\rm G})$. Conversely, one can recover
${\rm G}_x$ from $\vartheta(x)$ as its \emph{holomorphic envelope}
\cite[Proposition 2.4 (ii)]{RTW1}, which is to say:
$${\rm G}_x = \{y \in {\rm G}^{\rm an} \ ; \ \forall f \in
\mathcal{O}({\rm G}), \ |f(y)| \leqslant |f(\vartheta(x))|\}.$$
This can be phrased equivalently in terms of multiplicative norms on $\mathcal{O}({\rm G})$ by saying that one recovers the affinoid algebra of ${\rm G}_x$ as the completion of the normed $k$-algebra $\left(\mathcal{O}({\rm G}), |.|(\vartheta(x) \right)$.

\smallskip 

Let us say that a $k$-affinoid subgroup ${\rm H}$ of ${\rm G}^{\rm an}$ is of
\emph{Chevalley type} (or a \emph{Chevalley $k$-affinoid subgroup}) if there exists a $k^\circ$-Chevalley semi-simple group $\mathcal{H}$ \cite[Expos\'e XXIII \S 5]{SGA3} with $\mathcal{H} \otimes_{k^\circ} k \simeq G$ and such that $H$ is the generic fibre of the formal completion of $\mathcal{H}$ along its special fibre.  More generally, we will say that an affinoid subgroup ${\rm H}$ of ${\rm G}^{\rm an}$ is
\emph{potentially} of Chevalley type if there exists an affinoid extension $K/k$ such that ${\rm H}
\widehat{\otimes}_k K$ is a Chevalley affinoid subgroup of ${\rm G}^{\rm an} \widehat{\otimes}_k K$ ; by an
\emph{affinoid} extension, we simply mean that $K$ is a
non-Archimedean field which is a $k$-affinoid algebra (see \cite[Appendix A]{RTW1}; this restriction allows to recover $k$-affinoid algebras from $K$-affinoid algebras equipped with a descent datum). By
construction, the $k$-affinoid subgroup ${\rm G}_x$ attached to a
point $x$ of $\mathcal{B}({\rm G},k)$ is always potentially of
Chevalley type.

\smallskip 

For a point $z$ of ${\rm G}^{\rm an}$, let us define its
\emph{holomorphic envelope} by $${\rm G}(z) = \{y \in
     {\rm G}^{\rm an} \ ; \ \forall f \in \mathcal{O}({\rm G}),
     \ |f(y)| \leqslant |f(z)|\}.$$
The above discussion brings out a first condition fulfilled by
any point of ${\rm G}^{\rm an}$  belonging to the image of
$\vartheta$.

\medskip
 
\noindent
\textsc{First condition.---}~ 
\emph{The holomorphic envelope of $z$ is a $k$-affinoid subgroup potentially of Chevalley type.}

\medskip

It is easily checked that a point satisfying this condition does appear in the image
of $\vartheta$ over some non-Archimedean extension of $k$, as the next Lemma shows. 

\begin{Lemma}
Let $z$ be a point of ${\rm G}^{\rm an}$ whose holomorphic envelope is
a $k$-affinoid subgroup potentially of Chevalley type. 

\vskip1mm
\begin{itemize}
 \item[(i)] For every non-Archimedean extension $K/k$, the Shilov boundary of the $K$-affinoid domain $G(z) \widehat{\otimes}_k K$ is a singleton $\{z_K\}$.
 \vskip1mm
 \item[(ii)] There exists a non-Archimedean extension $K/k$ such that the point $z_K$ belongs to the image of $\vartheta_K$.
\end{itemize}
\end{Lemma}

\smallskip 


\medskip 

\noindent\emph{Proof.---}~ (i) Let $K/k$ be a non-Archimedean extension. By assumption, there exists an affinoid extension $K_0/k$ such that ${\rm G}(z) \widehat{\otimes}_k K_0$ is a $K_0$-affinoid subgroup of Chevalley type in ${\rm G}^{\rm an} \widehat{\otimes}_k K_0$. We consider a non-Archimedean extension $L/k$ containing both $K_0$ and $K$. The set $${\rm G}(z) \widehat{\otimes}_k L = ({\rm G}(z) \widehat{\otimes}_k K_0) \widehat{\otimes}_{K_0} L$$
is a $L$-affinoid subgroup of Chevalley type in ${\rm G}^{\rm an} \widehat{\otimes}_k L$, hence is the generic fibre of a formal scheme with geometrically integral special fibre; it implies that the Shilov boundary of ${\rm G}(z) \widehat{\otimes}_k L$ is a singleton $\{z_L\}$. Since the canonical projection map sends the Shilov boundary of ${\rm G}(z) \widehat{\otimes}_k L$ onto the Shilov boundary of ${\rm G}(z) \widehat{\otimes}_k K$, the latter is also a singleton.
\vskip1mm
(ii) We consider again an affinoid extension $K/k$ such that ${\rm G}(z) \widehat{\otimes}_k K$ is a
$K$-affinoid subgroup of Chevalley type in ${\rm G}^{\rm an}
\widehat{\otimes}_k K$. 

\smallskip 

The $K$-affinoid Chevalley subgroup ${\rm G}(z)
\widehat{\otimes}_k K$ is the stabilizer of a unique point $x$
of $\mathcal{B}({\rm G},K)$, hence $({\rm G}_K)_x = {\rm
  G}(z) \widehat{\otimes}_k K$ and therefore $\vartheta_{K}(x)
= z_{K}$. We used the fact
that any $K$-affinoid Chevalley subgroup of ${\rm G}^{\rm an} \widehat{\otimes}_k K$
occurs as the stabilizer of some point in the building. To see this, just pick
a special vertex; its stabilizer is a $K$-affinoid Chevalley subgroup of
${\rm G}^{\rm an} \widehat{\otimes}_k K$, and any two of them are
${\rm G}(K)$-conjugate.
\hfill $\square$

\section{Galois-fixed points in buildings}
\label{ss-Galois}

It is clear that the above condition does not suffice to
characterize the image of $\vartheta$. Indeed, consider a finite Galois
extension $k'/k$ and pick a point $x'$ in $\mathcal{B}({\rm G},k')$
which is fixed under the natural action of ${\rm Gal}(k'|k)$ on the
building. Let $z$ denote the image of $\vartheta_{k'}(x')$ under the
canonical projection of ${\rm G}_{k'}^{\rm an}$ onto ${\rm G}^{\rm
  an}$. The $k'$-affinoid subgroup $({\rm G}_{k'})_{x'}$ is equipped with a
Galois descent datum, from which one deduces: $({\rm G}_{k'})_{x'} =
{\rm G}(z) \otimes_k k'$. It follows that ${\rm G}(z)$ is a
$k$-affinoid subgroup potentially of Chevalley type. Now, if the field
extension $k'/k$ is wildly ramified, then the inclusion of
$\mathcal{B}({\rm G},k)$ into the set of Galois-fixed points in
$\mathcal{B}({\rm G},k')$ is strict in general (see example below); therefore, if we pick a Galois-fixed point $x'$ outside $\mathcal{B}({\rm
  G},k)$, then $z$ does not belong to the image of $\vartheta$. 
  Therefore we see, up to illustrating the strict inclusions mentioned above, that this argument shows that we need additional conditions to describe the image of $\vartheta$ in ${\rm G}^{\rm an}$. 
  
\smallskip 

We want to illustrate this discussion by looking at an elementary
example. Let us consider the group ${\rm G} = {\rm SL}_2$ over
  some discretely valued field $k$ and pick a finite Galois extension
  $k'$ of $k$. Via its canonical embedding in $\mathbf{P}^{1,{\rm an}}_{k'}$, the building $\mathcal{B}({\rm G},k')$ can be identified with the convex hull of $\mathbf{P}^1 (k')$ inside
$\mathbf{P}^{1, {\rm an}}_{k'}$ with $\mathbf{P}^1 (k')$ omitted, {\it i.e.} with the subset $$\bigcup_{a \in k}
  \eta_a\left(\mathbf{R}_{>0}\right),$$ where $\eta_a$ denotes the map
  from $\mathbf{R}_{>0}$ to $\mathbf{A}^{1,{\rm an}}_k$ sending $r$ to the maximal point
  of the ball of radius $r$ centered in $a$. The Galois
  action on $\mathcal{B}({\rm G},k')$ is induced by the Galois
  action on $\mathbf{P}^{1,{\rm an}}_{k'}$, and the sub-building $\mathcal{B}({\rm G},k)$
  is the image of paths $\eta_a$ with $a \in k$. Since the field
  $k$ is discretely valued -- hence spherically complete -- there
  exists a well-defined Galois-equivariant retraction
  $$ \tau : \mathbf{P}^{1, {\rm an}}_{k'} - \mathbf{P}^{1,{\rm
      an}}(k) \longrightarrow \mathcal{B}({\rm G},k)$$ defined
  by sending a point $x$ to the maximal point of the smallest ball
  with center in $k$ containing $x$.

\smallskip 

Using this picture, one easily sees how a Galois-fixed point can appear
  in $\mathcal{B}({\rm G},k') -\mathcal{B}({\rm G},k)$. It suffices to find an element
  $\alpha$ of $k'$ such that all the paths
  $\eta_{\alpha^g}(\mathbf{R}_{>0})$ issued from conjugates $\alpha^g$ of $\alpha$ intersect at some point
  distinct form $\tau(\alpha)$; since the Galois action permutes these paths, their meeting point $x'$ will be fixed. Note that we have $$x' =
  \eta_{\alpha}(r) \ \ \ {\rm and} \ \ \ \tau(\alpha) = \eta_{\alpha}(r'),$$ where $r =
  \max \{|\alpha^g - \alpha| \ ; \ g \in {\rm
    Gal}(k'|k)\}$ is the diameter of the Galois orbit of $\alpha$ and $r' = \min \{|\alpha -
  a| \ ; \ a \in k\}$ is the distance from $\alpha$ to $k$.

\smallskip  

Let $k'$ be any totally ramified finite Galois extension of $k$. It is
well-known that $k'$ can be realized as the splitting field of some
Eisenstein polynomial ${\rm P}({\rm T})={\rm T}^e + a_{e-1} {\rm T}^{e-1} + \ldots +
a_1 {\rm T} + a_0$, where $|a_i| \leqslant |a_0| < 1$ for all
$i$ and $|k^{\times}| = |a_0|^{\mathbf{Z}}$ (see for example Chapter 7, Theorem 7.1, in \cite{Cassels}). The group
$|{k'}^{\times}|$ is generated by $|\alpha| =
|a_0|^{1/e}$ for any root $\alpha$ of ${\rm P}$.
We have ${\rm d}(\alpha,k) = |\alpha|$ and all conjugates of $\alpha$ are
contained in the closed ball ${\rm E}(\alpha, |\alpha|) = {\rm E}(0,|\alpha|)$. The
endomorphism of $\mathbf{A}^{1,{\rm an}}_{k'}$ defined by ${\rm P}({\rm
  T})$ maps this ball onto the closed ball ${\rm E}(0,
|a_0|)$. In order to study the induced map ${\rm E}(0,|\alpha|) \rightarrow {\rm E}(0,|a_0|)$, set ${\rm U}
= {\rm T}/\alpha$ and write $${\rm Q}({\rm U}) = \frac{1}{a_0} {\rm
  P}(\alpha{\rm U}) = \frac{\alpha^e}{a_0} {\rm U}^e + \frac{a_{e-1} \alpha^{e-1}}{a_0} {\rm U}^{e-1} +
\ldots + \frac{a_1}{a_0} \alpha {\rm U} + 1.$$ Since
$|a_i|.|\alpha|^{i}< |a_i| \leqslant |a_0|$ for any $i \in \{1,
\ldots, e-1\}$, the polynomial ${\rm Q}$ reduces to
$\widetilde{\frac{\alpha^e}{a_0}} {\rm U}^e+1 = 1 - {\rm U}^e$ in
$\widetilde{k'}[{\rm U}]$. It follows that the following four
conditions are equivalent:

\begin{itemize}
\item[-] all paths $\eta_{\alpha^g}(\mathbf{R}_{>0})$, for
  $g \in {\rm Gal}(k'|k)$, intersect outside $\mathcal{B}({\rm G},k)$;
\item[-] all roots of ${\rm P}$ are contained in the \emph{open}
  ball ${\rm D}(\alpha, |\alpha|)$;
\item[-] all roots of ${\rm Q}$ are contained in the \emph{open} ball
  ${\rm D}(1,1)$;
\item[-] $e$ vanishes in $\widetilde{k}$.
\end{itemize}

\begin{figure}[h]
\centering
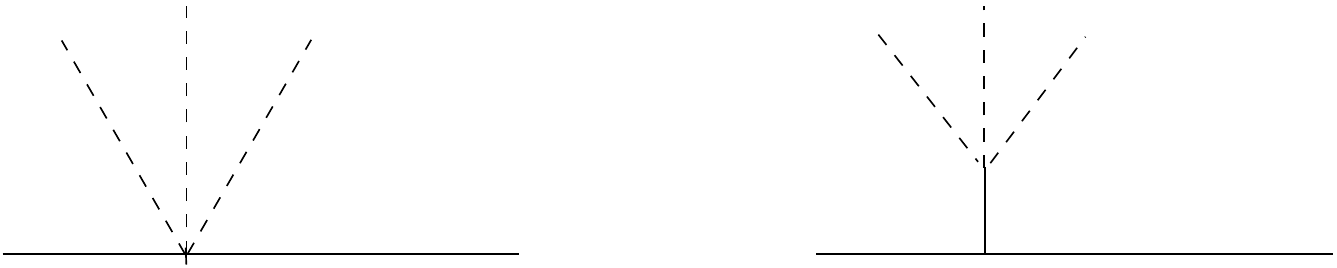
\end{figure}

In particular, for any totally ramified (finite) Galois extension $k'/k$, the building
$\mathcal{B}({\rm G},k)$ is strictly smaller than the
set of Galois-fixed points in $\mathcal{B}({\rm G},k')$ if and only if
$[k':k]$ is divisible by the residue characteristic.
We refer to \cite{PrasadFP} for a recent contribution to the question of describing fixed-point sets of finite group actions in Bruhat-Tits theory. 

\begin{Example} Let $k=\mathbf{Q}_2$ and $k' =
  \mathbf{Q}_2(\alpha)$, where $\alpha^2=2$. The two paths
  $\eta_{\alpha}(\mathbf{R}_{>0})$ and
  $\eta_{-\alpha}(\mathbf{R}_{>0})$ intersect $\mathcal{B}({\rm G},k)$
  along the image of $[2^{-1/2},\infty)$, whereas they meet along the image of
  $[2^{-3/2},\infty)$. The whole interval
      $\eta_\alpha\left([2^{-3/2},2^{-1/2})\right)$ consists of Galois-fixed
        points lying outside $\mathcal{B}({\rm G},k)$. In general, Rousseau gave an upper bound for the distance of a Galois-fixed point in $\mathcal{B}({\rm G},k')$ to $\mathcal{B}({\rm G},k)$ in terms of the ramification of $k'/k$ \cite[Proposition 5.2.7]{RousseauOrsay}.
\end{Example}

\section{Apartments}
\label{ss-Recovering apartments}

As we have just seen, the characterization of the image of the Bruhat-Tits building of $\rm G$ inside ${\rm G}^{\rm an}$ requires an additional condition. 
The goal of this section is to formulate such a condition which involves maximal tori of ${\rm G}$. 
We will have to make use of the following fact.

\begin{Lemma} \label{tori} Let ${\rm T}$ be a torus over $k$.
\begin{itemize}
\item[(i)] Its analytification ${\rm T}^{\rm an}$ contains a largest bounded subgroup ${\rm T}^1$. This is an affinoid subgroup, which coincides with the affinoid domain cut out by the equations $|\chi|=1$ for all $\chi \in {\rm X}^*({\rm T})$ if ${\rm T}$ is split.
\item[(ii)] The Shilov boundary of ${\rm T}^1$ is reduced to a single point.
\end{itemize}
\end{Lemma}

\smallskip 

We will write ${\rm o}_{\rm T}$ for the unique Shilov boundary point of ${\rm T}^1$.

\medskip 

\noindent\emph{Proof.---}~
We consider first the case of a split torus. If $\chi_1, \ldots, \chi_n$ is a basis of characters of ${\rm T}$, then the equations $|\chi_1|=1,\ldots, |\chi_n|=1$ cut out a $k$-affinoid subgroup ${\rm T}^1$ of ${\rm T}^{\rm an}$ over which $|\chi|=1$ for any character $\chi \in {\rm X}^*({\rm T})$. Let $K/k$ be a non-Archimedean extension and $\Gamma$ a bounded subgroup of ${\rm T}(K)$. For any character $\chi$ of ${\rm T}$, both $\chi(\Gamma)$ and $(-\chi)(\Gamma) = \chi(\Gamma)^{-1}$ are bounded subgroups of $K^{\times}$, hence $|\chi(\Gamma)|=1$ and $\Gamma \subset {\rm T}^1({\rm K})$. The Shilov boundary of ${\rm T}^1$
 is reduced to a point since the reduction of ${\rm T}^1$, a torus over $\widetilde{k}$, is irreducible.

\smallskip 

In general, pick a finite Galois extension $k'/k$ splitting ${\rm T}$ and set ${\rm T}_{k'} = {\rm T} \otimes_k k'$. The affinoid subgroup ${\rm T}_{k'}^1$ of ${\rm T}_{k'}^{\rm an}$ is stable under the natural Galois action on ${\rm T}_{k'}^{\rm an}$, hence its descends to a $k$-affinoid subgroup ${\rm T}^1$ of ${\rm T}$ such that ${\rm T}_{k'}^1 = ({\rm T}^1) \otimes_k k'$. Finally, since the Shilov boundary of ${\rm T}^1 \otimes_k k'$ is the preimage of the Shilov boundary of ${\rm T}^1$ under the canonical projection, we see that ${\rm T}^1$ contains a unique Shilov boundary point.
\hfill$\Box$

\smallskip 

Let us now go back to our discussion toward a characterization of the building inside ${\rm G}^{\rm an}$. 
We assume temporarily that the group ${\rm G}$ is split. 
Given a point $z$ in ${\rm G}^{\rm an}$ whose holomorphic envelope, again denoted by ${\rm G}(z)$, is a $k$-affinoid subgroup potentially of Chevalley type, let us consider a non-Archimedean extension $K/k$ such that the canonical
lifting $z_K$ of $z$ belongs to the image of $\vartheta_K$ and denote by $x$ its preimage:
$z_K = \vartheta_K(x)$. Since the group ${\rm G}$ is split, the embedding $\mathcal{B}({\rm G},k) \hookrightarrow \mathcal{B}({\rm G},K)$ identifies the left-hand side with the union of apartments of all maximal split tori in ${\rm G}_K$ which are defined over $k$. 
Therefore, in order to guarantee that the point $z$ itself belongs to the  image of $\vartheta$, we should require that $x$ belongs to the apartment of a maximal split torus defined over $k$. 
The next proposition translates this additional condition into appropriate terms.

\begin{Prop}
\label{prop-bounded_torus} 
Let ${\rm S}$ be a maximal split torus and let $x$ be a point of $\mathcal{B}({\rm G},k)$. 
The following conditions are equivalent:
\begin{itemize}
\item[(i)] the point $x$ belongs to the apartment ${\rm A}({\rm S},k)$;
\item[(ii)] for every non-Archimedean extension $K/k$, the point $x$ is fixed by the ${\rm S}^1(K)$-action on $\mathcal{B}({\rm G},K)$;
\item[(iii)] the affinoid subgroup ${\rm G}_x$ of ${\rm G}^{\rm an}$ contains ${\rm S}^1$.
\end{itemize}
\end{Prop}

Note that the proof will provide another equivalent statement, namely: 

\begin{itemize}
\item[(ii)$'$] \emph{for every finite unramified extension $k'/k$, the point $x$ is fixed by the ${\rm S}^1(k')$-action on $\mathcal{B}({\rm G},k')$}. 
\end{itemize}

\medskip

\noindent 
{\it Proof.---}~
Equivalence of assertions (ii) and (iii) follows immediately from the definition of the affinoid group ${\rm G}_x$.

\smallskip 

For any non-Archimedean extension $K/k$, the action of ${\rm S}(K)$ on $\mathcal{B}({\rm G},K)$ preserves the apartment ${\rm A}({\rm S},K)$ and the induced action of the maximal bounded subgroup ${\rm S}^1(K)$ is trivial  \cite[Proposition 4.18 (iii)]{RTW3}. 
Hence (i) implies (ii).

\smallskip 

The only non trivial point to prove is that (ii) implies (i). 
We argue by contraposition: given a point $x$ of the building which does not belong to the apartment ${\rm A}={\rm A}({\rm S},k)$, we will exhibit a finite unramified extension $k'/k$ and an element $s$ of ${\rm S}^1(k')$ such that $s \cdot x \neq x$ in $\mathcal{B}({\rm G},k')$.

Recall that the building $\mathcal{B}({\rm G},k)$ is also a complete metric space satisfying a non-positive curvature property \cite[2.5 and 3.2]{BT1a}, called by now the \emph{${\rm CAT}(0)$-property} (see \cite[Part II]{BriHae}). 
Since the apartment ${\rm A}$ is a closed and convex subset, it contains a unique point $x'$ with $${\rm d}(x,{\rm A}) = {\rm d}(x,x').$$ 
Moreover, there exists a unique geodesic segment $[x,x']$ between those two points and $$[x,x'] \cap {\rm A} = \{x'\}.$$  Consider any apartment ${\rm A}'$ containing both $x$ and $x'$; it contains $[x,x']$ and intersects ${\rm A}$ along a closed convex subset ${\rm C}$. The point $x'$ coincides with the projection in ${\rm A}'$ of $x$ to ${\rm C}$, hence it lies in the boundary of ${\rm C}$ since $x \not\in {\rm A}$. Finally, since ${\rm C}$ is the intersection of a finite number of half-apartments \cite[Proposition 9.1]{RousseauGrenoble}, we conclude that $x'$ belongs to the wall ${\rm H}$ in ${\rm A}$ defined by the vanishing of some affine root $a$. Let ${\rm H}^+$ denote the half-apartment on which $a$ is non-negative and write $a=\alpha+\lambda$, where
the linear part $\alpha$ belongs to the root system $R({\rm
  T},{\rm G})$ and $\lambda$ is an element of ${\rm log}|k^\times|$.

\smallskip 

Recall now that the unipotent root group ${\rm U}_\alpha$ is endowed with a separated, exhaustive and decreasing filtration $\{{\rm U}_{\alpha,s}\}_{s \in \mathbf{R}}$ by affinoid subgroups \cite[4.2.1]{RTW3}. 
The subgroup ${\rm U}_a = {\rm U}_{\alpha,\lambda}$ corresponding to the affine root $a = \alpha + \lambda$ has the following geometric interpretation: for every non-Archimedean extension $K/k$, the action of ${\rm U}_{a}(K)$ on $\mathcal{B}({\rm G},K)$ fixes pointwise the half-apartment ${\rm H}^+$ and is transitive on the set of apartments which intersect ${\rm A}$ along ${\rm H}^+$. 
Moreover, if we set 
$${\rm U}_{\alpha, \lambda}^+ = \bigcup_{r>\lambda} {\rm U}_{\alpha,r},$$
then ${\rm U}_{\alpha, \lambda}(K)/{\rm U}_{\alpha,\lambda}^+(K)$ is in bijection with equivalence classes of apartments in $\mathcal{B}({\rm G},K)$ containing ${\rm  H}^+$, where two such apartment are said to be
equivalent they intersect along a neighborhood of ${\rm H}^+$. 

\smallskip

The torus ${\rm S}$ acts on ${\rm U}_{\alpha}$ by conjugation. The
bounded subgroup ${\rm S}^1$ preserves each step of the filtration and
there is a non-canonical bijection $${\rm
  U}_{\alpha,\lambda}(K)/{\rm U}_{\alpha,\lambda}^+(K)
\simeq \widetilde{K}$$ such that the action of an element $s \in {\rm S}^1(K)$ on
the left-hand side corresponds to multiplication by
$\widetilde{\alpha(s)} \in \widetilde{K}^{\times}$ on the
right-hand side. Note that this condition implies that the unit
element on the left-hand side corresponds to $0$ on the right-hand side.

\smallskip 

To conclude the proof, observe that there exists a unique element $\xi \in \widetilde{k}$ satisfying the following property: for every $u \in {\rm U}_{\alpha,\lambda}(k)$ whose class modulo ${\rm U}^+_{\alpha,\lambda}(k)$ corresponds to $\xi$, the apartment $u \cdot {\rm A}$ contains a germ of $[x,x']$ at $x'$, {\it i.e.} we may write $u \cdot {\rm A} \cap [x,x'] = [x,y]$ for some $y \in [x,x')$. We
have $\xi \neq 0$ since $x$ does not belong to ${\rm A}$. Now, if
$x$ is fixed by ${\rm S}^1(k)$, then $\widetilde{\alpha(s)} \cdot
\xi = \xi$ and therefore $\widetilde{\alpha(s)}=1$ for any $s \in
{\rm S}^1(k)$. This means that the character $\widetilde{\alpha}$ of
the $\widetilde{k}$-torus $\widetilde{\rm S^1}$ is trivial at the
level of $\widetilde{k}$-points. Since the character $\widetilde{\alpha}$ is non-trivial,
there exists a finite separable extension $\ell$ of $\widetilde{k}$
such that $\widetilde{\alpha} \neq 1$ on $\widetilde{\rm
  S}(\ell)$. Denoting by $k'$ the corresponding unramified
extension of $k$, we conclude that $x$ is not fixed by the whole group ${\rm S}^1(k')$.
\hfill$\square$

\begin{Remark} 
The arguments above are close to the ones proving that the apartment ${\rm A}({\rm S},k)$ coincides with the fixed-point set of ${\rm S}^1(k)$ if $\widetilde{k}$ contains at least $4$ elements \cite[5.1.37]{BT1b}. 
\end{Remark}

In the split case, the discussion above shows precisely which
additional condition is required in order to characterize the image
of $\vartheta$ in ${\rm G}^{\rm an}$: \emph{there exists a maximal split torus ${\rm S}$
in ${\rm G}$ such that ${\rm G}(z) \cap {\rm S} = {\rm S}^1$}.
In general, we impose the following:

\medskip
 
\noindent
\textsc{Second condition.---}~ 
\emph{There exists a maximal torus ${\rm T}$ in ${\rm G}$ which contains a maximal split torus and such that ${\rm G}(z) \cap {\rm T} = {\rm T}^1$.}

\medskip

Once this second condition has been introduced naturally, we are in good position to characterize the image of $\vartheta$ in the analytic space ${\rm G}^{\rm an}$ of ${\rm G}$.

\begin{Thm}
\label{Thm-image}
The image of the canonical embedding $\vartheta : \mathcal{B}({\rm
  G},k) \hookrightarrow {\rm G}^{\rm an}$ is the subset of points $z$
satisfying the following two conditions:
\begin{itemize}
\item[1.] the holomorphic envelope ${\rm G}(z)$ of $z$ is a
  $k$-affinoid subgroup potentially of Chevalley type;
\item[2.] there exists a maximal torus ${\rm T}$ of ${\rm G}$
  containing a maximal split torus and such that ${\rm G}(z) \cap {\rm T}^{\rm an}$ is the maximal affinoid subgroup
  ${\rm T}^1$ of ${\rm T}^{\rm an}$.
\end{itemize}
\end{Thm}

\smallskip 

\noindent\emph{Proof.---}~ 
We have already seen that the first condition is necessary. 
The same holds for the second one. Given a point $x \in \mathcal{B} ({\rm G}, k)$, there exists a maximal split torus ${\rm S}$ and a maximal torus ${\rm T}$ containing ${\rm S}$ such that $x \in {\rm A}({\rm S}, k) \subset {\rm A}({\rm T},
k')$, where $k'$ is a finite extension of $k$ which splits ${\rm T}$. 
It follows that ${\rm G}_x \otimes_k k'$ contains the bounded torus ${\rm T}^1 \otimes_k k'$, hence ${\rm T}^1 \subset {\rm G}_x$.

\smallskip 

Now, let us show that the two conditions are sufficient. We consider a
point $z \in {\rm G}^{\rm an}$ satisfying the two
conditions and pick a finite Galois extension $k'/k$ spliting the maximal torus $T$ given by the second condition, as well as a non-Archimedean extension $K / k'$ and a point $w \in
\mathcal{B}({\rm G}, K)$ such that $z_K = \vartheta_K (w)$, where $z_K$
denotes the (unique) Shilov boundary point of $G(w) \widehat{\otimes}_k K$. 

Since
$${\rm T}^{\rm an}_K \cap \left( {\rm G}_K \right)_x =
   {\rm T}^{\rm an}_K \cap {\rm G} (z)_K = \left(
   {\rm T}^{\rm an} \cap {\rm G} (z) \right)_K = \left(
   {\rm T}^1 \right)_K = \left( {\rm T}_K \right)^1,$$
it follows from Proposition \ref{prop-bounded_torus} that $w$ belongs to the
apartment of ${\rm T}_K$. The canonical embedding of $\mathcal{B}({\rm G},k')$ into $\mathcal{B}({\rm G},K)$ induces a bijection between ${\rm A}({\rm T},k')$  and ${\rm A}({\rm T},K)$, hence $w$ comes from a point $y$ of $\mathcal{B}({\rm G},k')$ contained in ${\rm A}({\rm T},k')$. Moreover $({\rm G}_{k'})_y = H_{k'}$, for both sides coincide with $({\rm G}_{K})_w$ after base change to $K$. In particular, this shows that $y$ is fixed by ${\rm Gal}(k'|k)$. Since ${\rm A}({\rm T},k')^{{\rm Gal}(k'|k)}$ is the image of ${\rm A}(S,k)$ in $\mathcal{B}({\rm G},k')$, we get a point $x$ in ${\rm A}({\rm S},k)$ such that $$({\rm G}_{z_{k'}}) \otimes_k k' = ({\rm G}_{k'})_y = H_{k'},$$ hence $H = {\rm G}_x$ and $z = \vartheta(x)$. \hfill$\Box$

\section{A reformulation in terms of norms}
\label{ss-Reformulation norms}

The above characterization of points of ${\rm G}^{\rm an}$ lying inside the building $\mathcal{B}({\rm G},k)$ (identified with its image by the canonical map $\vartheta$) can be conveniently rephrased in terms of (multiplicative) $k$-norms on the coordinate algebra $\mathcal{O}({\rm G})$. As we are going to explain, we need here to make the additional assumption that the group ${\rm G}$ splits over a \emph{tamely ramified} extension of $k$.

\smallskip 
Let us start by recalling the notion of a \emph{universal point} \footnote{This notion was introduced by Berkovich, who used the adjective \emph{peaked} \cite[5.2]{Ber1}; its study was carried on by Poineau, who preferred the adjective \emph{universal} \cite{Poineau}.}. Let $z$ be a point in ${\rm G}^{\rm an}$, seen as a multiplicative $k$-seminorm on $\mathcal{O}({\rm G})$. For a given non-Archimedean field extension $K/k$, there is a natural $K$-seminorm $||.|| = x \otimes 1$ on $\mathcal{O}({\rm G}) \otimes_k K$, defined by $$||a|| = \inf \max_i |a_i(z)|\cdot |\lambda_i|$$
where the infimum is taken over the set of all expressions $\sum_i a_i \otimes \lambda_i$ representing $a$, with $a_i \in \mathcal{O}({\rm G})$ and $\lambda_i \in K$. The point $z$ is said to be \emph{universal} if, for any non-Archimedean field extension $K/k$, the above $K$-seminorm on $\mathcal{O}({\rm G}) \otimes_k K$ is multiplicative. 
One then writes $z_K$ for the corresponding point in ${\rm G}^{\rm an} \widehat{\otimes}_k K$.


\begin{Remarks}
\label{universal}
\begin{enumerate}
\item[1.]~Obviously, points of ${\rm G}^{\rm an}$ coming from $k$-rational points of ${\rm G}$ are universal.
\item[2.]~Let $x \in {\rm G}^{\rm an}$ be universal. For any finite Galois extension $k'/k$, the canonical extension $x_{k'}$ of $x$ to ${\rm G}^{\rm an} \otimes_k k'$ is invariant under the action of ${\rm Gal}(k'/k)$: indeed, the $k'$-norm $x \otimes 1$ on $\mathcal{O}({\rm G}) \otimes_k k'$ is Galois invariant.
\item[3.]~If ${\rm G}$ is split, then $\vartheta(x)$ is universal for every point $x$ of $\mathcal{B}({\rm G},k)$, as follows readily from the explicit description of  $\vartheta(x)$ as a norm on the coordinate algebra of a big cell \cite[Proposition 2.6]{RTW1}. More generaly, the same is true if ${\rm G}$ splits over a tamely ramified extension $k'/k$ by \cite[Lemma A.10]{RTW1} and \cite[Erratum]{RTW3}.
\item[4.]~Given a torus ${\rm T}$ over $k$, we denote by $\sigma_{{\rm T}}$ the (unique) Shilov boundary point of the maximal affinoid subtorus ${\rm T}^1$ of ${\rm T}^{\rm an}$. This point is universal if $T$ splits over a tamely ramified extension of $k$, but not in general as checked by Mayeux  \cite[Proposition 10.2]{Mayeux}.
\end{enumerate}
\end{Remarks}

\vskip2mm \noindent \textsc{Setting} --- \emph{As already mentioned, the semisimple group ${\rm G}$ considered in this last section is assumed to split over a tamely ramified extension of $k$, so that all points in the image of $\vartheta$ are universal.}
\vskip2mm
Let ${\rm G}^{\rm an}_u$ denote the subset of universal points in ${\rm G}^{\rm an}$. 
Following Berkovich \cite[5.2]{Ber1}, there is a natural \emph{monoid} structure on ${\rm G}^{\rm an}_u$ extending the group structure on ${\rm G}(k)$. Given any two points $g, h \in {\rm G}_u^{\rm an}$, the seminorm $g \otimes h$ on $\mathcal{O}({\rm G}) \otimes_k \mathcal{O}({\rm G})$ is multiplicative and one defines $g \ast h$ as the induced multiplicative seminorm on $\mathcal{O}({\rm G})$ via the comultiplication map $\Delta : \mathcal{O}({\rm G}) \rightarrow \mathcal{O}({\rm G}) \otimes_k \mathcal{O}({\rm G})$. This binary operation is associative, with unit the element $1 \in {\rm G}(k)$; moreover, it is obvious that we recover the group law if $g$ and $h$ belong to ${\rm G}(k)$.

\smallskip

More generally, given an (analytic) action of ${\rm G}^{\rm an}$ on some $k$-analytic space ${\rm X}$, one can define in a similar way an action of the monoid ${\rm G}^{\rm an}_u$ on the topological space underlying ${\rm X}$, which extends the action of ${\rm G}(k)$. In particular, $x \ast y$ is well-defined for any two points $x,y$ in ${\rm G}^{\rm an}$ such that one of them is universal.

\smallskip 

Finally, we also recall that one can define a partial order on the set underlying ${\rm G}^{\rm an}$ as follows:
$$x \preccurlyeq y \ \ \ {\rm if \ and \ only \ if} \ \ \forall a \in \mathcal{O}({\rm G}), \ |a(x)| \leqslant |a(y)|.$$

\smallskip

We can now give the following description of the building inside ${\rm G}^{\rm an}$.

\begin{Thm}
\label{characterization2} 
Let $x \in {\rm G}^{\rm an}$. 
Then the point $x$ belongs to the image of the canonical embedding $\vartheta : \mathcal{B}({\rm G},k) \hookrightarrow {\rm G}^{\rm an}$ if, and only if, the following conditions are satisfied: 
\begin{itemize}
\item[(i)] $x$ is universal;
\item[(ii)] $x \ast x \preccurlyeq x$ and ${\rm inv}(x) \preccurlyeq x$;
\item[(iii)] there exists a maximal torus ${\rm T}$ containing a maximal split torus such that ${\rm o}_{\rm T} \preccurlyeq x$;
\item[(iv)] $x$ is maximal with respect to the three conditions above.
\end{itemize}
\end{Thm}

\begin{Remarks}
\begin{enumerate}
\item[1.]~Conditions (i)-(iv) imply that $x$ is a multiplicative $k$-\emph{norm} on $\mathcal{O}({\rm G})$.
\item[2.]~One way to understand condition (ii) is to say that $x$ defines a $k$-multiplicative norm on the commutative Hopf algebra $\mathcal{O}({\rm G})$ with respect to which commultiplication and the antipode are bounded. 
One should observe that (ii) obviously implies $1_{\rm G} \preccurlyeq x$, so that the counit is also bounded.
\item[3.]~By Lemma \ref{*-order} below, the condition $1_{\rm G} \preccurlyeq x$ (resp. ${\rm inv}(x) \preccurlyeq x$) implies $x = 1_{\rm G} \ast x \preccurlyeq x \ast x$ (resp. $x \preccurlyeq {\rm inv}(x)$), hence condition (ii) could be replaced by $x \ast x = x$ and ${\rm inv}(x)=x$.
\end{enumerate}
\end{Remarks}

\begin{Lemma} 
\label{*-order} 
Let $x,y,x',y'$ be four points in ${\rm G}^{\rm an}$ such that both sets $\{x,y\}$ and $\{x',y'\}$ contain at least one universal point (so that $x \ast y$ and $x' \ast y'$ are well-defined). If $x \preccurlyeq x'$ and $y \preccurlyeq y'$, then $x \ast y \preccurlyeq x' \ast y'$ and ${\rm inv}(x) \preccurlyeq {\rm inv}(x')$.
\end{Lemma}

\smallskip 

\noindent
\emph{Proof.---}~
This follows directly from the formulas $$\forall a \in \mathcal{O}({\rm G}), \ \ \ |a(x \ast y)| = \inf \max_i |b_i(x)| \cdot |c_i(y)| \ \ {\rm and } \ \ |a(x' \ast y')| = \inf \max_i |b_i(x')| \cdot |c_i(y')|,$$
where the infimum is taken over the set of all expressions $\sum_i b_i \otimes c_i$ representing $\Delta(a)$, and 

\medskip 

\centerline{\hfill $\forall a \in \mathcal{O}({\rm G}), \ \ \ |a({\rm inv}(x))| = |{\rm inv}(a)(x)| \leqslant |{\rm inv}(a)(x')| = |a({\rm inv}(x'))|$. \hfill $\Box$}

\medskip

\begin{Lemma} 
\label{*-group} 
Let $x$ be a universal point of ${\rm G}^{\rm an}$. We again denote by ${\rm G}(x) = \{z \in {\rm G}^{\rm an} \ ; \ z \preccurlyeq x\}$ its holomorphic envelope. 
Then the following conditions are equivalent:
\begin{itemize}
\item[(i)] The envelope ${\rm G}(x)$ is a subgroup object of ${\rm G}^{\rm an}$, {\it i.e. }${\rm G}(x)({\rm K})$ is a subgroup of ${\rm G}({\rm K})$ for any non-Archimedean extension $K/k$;
\item[(ii)] The point $x$ satisfies: $1_{\rm G} \preccurlyeq x, \ \ \ {\rm inv}(x) \preccurlyeq x \ \ {\rm and} \ \ \ x \ast x \preccurlyeq x$. 
\end{itemize}
Moreover, ${\rm G}(x)$ is bounded in ${\rm G}^{\rm an}$.
\end{Lemma}

\smallskip 

\noindent\emph{Proof.---}~
Assume that ${\rm G}(x)$ is a subgroup object of ${\rm G}^{\rm an}$. Since ${\rm G}(x)(k)$ is a subgroup of ${\rm G}(k)$, it contains the unit element $1_{\rm G}$ and therefore $1_{\rm G} \preccurlyeq x$. Now we consider the canonical point $\underline{x} \in {\rm G}(\mathcal{H}(x))$ lying over $x$; we have $|a(x)| = |a(\underline{x})|$ for any $a \in \mathcal{O}({\rm G})$, as well as $$|a(x \ast x)| = |a(\underline{x} \cdot \underline{x})| \ \ \ {\rm and} \ \ \ |a({\rm inv}(x))| = |a(\underline{x}^{-1})|.$$
Since $\underline{x} \cdot \underline{x}$ and $\underline{x}^{-1}$ belong to ${\rm G}(x)(\mathcal{H}(x))$, it follows that $$|a(x \ast x)| \leqslant |a(x)| \ \ \ {\rm and} \ \ \ |a({\rm inv}(x))| \leqslant |a(x)|$$
for all $a \in \mathcal{O}({\rm G})$, which is to say $x \ast x \preccurlyeq x$ and ${\rm inv}(x) \preccurlyeq x$.

\smallskip 

We assume now that $x$ is a universal point of ${\rm G}^{\rm an}$ satisfying the conditions of (ii). Obviously, ${\rm G}(x)$ contains the $k$-rational point $1_{\rm G}$. Given a non-Archimedean extension $K/k$ and elements $g,h \in {\rm G}(x)(K)$, we have for all $a \in \mathcal{O}({\rm G})$:
$$|a(g^{-1})| = |{\rm inv}(a)(g)| \leqslant |{\rm inv}(a)(x)| = |a({\rm inv}(x))|$$ and $$|a(gh)| = |\Delta(a)(g,h)| \leqslant \inf \max_i |b_i(g)|\cdot |c_i(h)| \leqslant \inf \max_i |b_i(x)|\cdot |c_i(x)| = |a(x \ast x)|,$$
where the infimum is taken over the set of all expressions $\sum_i b_i \otimes c_i$ representing $\Delta(a)$. Since ${\rm inv}(x) \preccurlyeq x$ and $x \ast x \preccurlyeq x$, we deduce $$|a(g^{-1})| \leqslant |a(x)| \ \ \ {\rm and} \ \ \ |a(gh)| \leqslant |a(x)|,$$ hence $g^{-1}, gh \in {\rm G}(x)({\rm K})$. This proves that ${\rm G}(x)$ is a subgroup object of ${\rm G}^{\rm an}$.

\smallskip 

Finally, boundedness is obvious: if $f_1, \ldots, f_n$ is a finite set generating $\mathcal{O}({\rm G})$ as a $k$-algebra, then we have $|f_i(y)| \leqslant \max_i |f_i(x)|$ for any point $y \in {\rm G}(x)$.
\hfill$\Box$

\medskip 

\noindent
\emph{Proof of Theorem \ref{characterization2}.---}~
In what follows, we identify the building $\mathcal{B}({\rm G},k)$ with its image in ${\rm G}^{\rm an}$ by the embedding $\vartheta$.

\smallskip 

If a point $x$ of ${\rm G}^{\rm an}$ belongs to $\mathcal{B}({\rm G},k)$, then $x$ is universal by Remark \ref{universal} and ${\rm G}(x)$ is a $k$-affinoid subgroup of ${\rm G}^{\rm an}$, hence $x \ast x \preccurlyeq$ and ${\rm inv}(x) \preccurlyeq x$ by Lemma \ref{*-group}. Moreover, there exists a maximal torus ${\rm T}$ containing a maximal split torus and such that ${\rm G}(x) \cap {\rm T}^{\rm an} = {\rm T}^1$, which amounts to saying that $x$ dominates the distinguished point ${\rm o}_{\rm T}$ of ${\rm T}$. Finally, consider a universal point $z \in {\rm G}^{\rm an}$ satisfying condition (ii) and dominating $x$ (which implies that $z$ dominates $1_{\rm G}$). For any non-Archimedean extension $K/k$, Lemma \ref{*-group} implies that ${\rm G}(z)(K)$ is a bounded subgroup of ${\rm G}(K)$ containing ${\rm G}(x)(K)$; by maximality of the latter, we deduce ${\rm G}(z)(K) = {\rm G}(x)(K)$, hence ${\rm G}(x) = {\rm G}(z)$ and $z=x$. We have thus checked that $x$ satisfies conditions (i)-(iv).

\smallskip 

Conversely, let $x$ be a point in ${\rm G}^{\rm an}$ satisfying conditions (i)-(iv). We observe that condition (iv) implies that $x$ dominates $1_{\rm G}$, hence it follows from Lemma \ref{*-group} that ${\rm G}(x)(K)$ is a bounded subgroup of ${\rm G}(K)$ for any non-Archimedean extension $K/k$. We are going to see that all these subgroups fix a common point in $\mathcal{B}({\rm G},k)$. By condition (iii), there exists a maximal torus ${\rm T}$ of ${\rm G}$ containing a maximal split torus ${\rm S}$ and such that $\sigma_T  \preccurlyeq x$, hence ${\rm T}^1(K) \subset {\rm G}(x)(K)$ for every non-Archimeden extension $K/k$.

\smallskip

We first assume that ${\rm T}$ is split and that $\widetilde{k}$ contains at least 4 elements. For any non-Archimedean field extension $K/k$, we identify ${\rm A}({\rm T},k)$ with ${\rm A}({\rm T},K)$ in $\mathcal{B}({\rm G},K)$, and we let $\mathcal{F}(K)$ denote the fixed-point set of ${\rm G}(x)(K)$ in $\mathcal{B}({\rm G},K)$. This is a non-empty closed subset, which is contained in the apartment ${\rm A}({\rm T},k)$ since ${\rm G}(x)(K)$ contains the group of units ${\rm T}^1(K)$ and there are at least 4 elements in $\widetilde{K}$. Considering in particular the extension $\mathcal{H}(x)/k$, we get a point $z \in {\rm A}({\rm T},k)$ fixed under the canonical element $\underline{x} \in {\rm G}(\mathcal{H}(x))$ lifting $x$. This means that $\underline{x}$ belongs to ${\rm G}(z)(\mathcal{H}(x))$, or equivalently that $x$ is contained in the $k$-affinoid subgroup ${\rm G}(z)$. This amounts to saying that $z$ dominates $x$, hence $x=z$ by maximality.

\smallskip 

In general, we consider a finite Galois extension $k'/k$ which splits ${\rm T}$ and such that $\widetilde{k'}$ contains at least four elements. It follows from the argument above that the canonical extension $x_{k'}$ of $x$ to ${\rm G}^{\rm an} \otimes_k k'$ belongs to the apartment of ${\rm T} \otimes_k k'$ in $\mathcal{B}({\rm G},k')$. Since this point is invariant under the action of ${\rm Gal}(k'/k)$ by Remark \ref{universal} 2, $x_{k'}$ belongs to the image of ${\rm A}({\rm S},k)$ in ${\rm A}({\rm T},k')$, and therefore $x$ belongs to ${\rm A}({\rm S},k)$. 
\hfill$\Box$

\vspace{1cm}

\begin{flushleft} \textit{Bertrand R\'emy} \\
Unit\'e de Math\'ematiques Pures et Appliqu\'ees (UMR 5669) \\
\'Ecole normale sup\'erieure de Lyon / CNRS / Inria \\
46 all\'ee d'Italie \\
F-69364 Lyon cedex 07 \\
\vspace{1pt}
bertrand.remy@ens-lyon.fr
\end{flushleft}

\vspace{0.1cm}
\begin{flushleft} \textit{Amaury Thuillier} \\
Institut Camille Jordan (UMR 5208) \\
Universit\'e de Lyon / Universit\'e Lyon 1 / CNRS \\
43 boulevard du 11 novembre 1918 \\
F-69622 Villeurbanne cedex \\
\vspace{1pt}
thuillier@math.univ-lyon1.fr
\end{flushleft}

\vspace{0.1cm}
\begin{flushleft}
\textit{Annette Werner} \\
Institut f\"ur Mathematik \\
Goethe-Universit\"at Frankfurt \\
Robert-Mayer-Str. 6-8 \\
D-60325 Frankfurt-am-Main \\
\vspace{1pt}
werner@math.uni-frankfurt.de
\end{flushleft}

\end{document}